\documentclass{article}

\usepackage{amsmath}
\usepackage{amsthm}
\usepackage{amssymb}
\usepackage{amscd}
\usepackage{ifthen}
\usepackage{epsfig,graphics,graphicx}
\usepackage{array,tabularx}
\usepackage{url}
\usepackage[usenames]{color}
\usepackage[pdftex,colorlinks=true,urlcolor=blue,linkcolor=blue,plainpages=false,pdfpagelabels]{hyperref}

\newtheorem{theorem}{Theorem}

\newtheorem{proposition}[theorem]{Proposition}

\newcommand{\bprop}{\begin{proposition}}
\newcommand{\eprop}{\end{proposition}}

\newcommand{\bcor}{\begin{corollary}}
\newcommand{\ecor}{\end{corollary}}

\newcommand{\be}{\begin{enumerate}}
\newcommand{\ee}{\end{enumerate}}
\newcommand{\bce}{\begin{center}}
\newcommand{\ece}{\end{center}}
\newcommand{\beq}{\begin{equation}}
\newcommand{\eeq}{\end{equation}}

\def\N{{\mathbb N}}
\newcommand{\ds}{\displaystyle}

\newcommand{\eps}{\varepsilon}
\newcommand{\limn}{\ds\lim_{n\to\infty}}

\newcounter{ct}

\newcommand{\details}[1]{}

\begin{document}

\title{A note on an alternative iterative method for nonexpansive mappings}
\author{Lauren\c{t}iu Leu\c{s}tean${}^{1,2}$,  Adriana Nicolae${}^{3,4}$ \\[0.2cm]
\footnotesize ${}^1$ Faculty of Mathematics and Computer Science, University of Bucharest,\\
\footnotesize Academiei 14,  P.O. Box 010014, Bucharest, Romania\\[0.1cm]
\footnotesize ${}^2$ Simion Stoilow Institute of Mathematics of the Romanian Academy,\\
\footnotesize P. O. Box 1-764, RO-014700 Bucharest, Romania\\[0.1cm]
\footnotesize ${}^3$ Department of Mathematics, Babe\c{s}-Bolyai University, \\
\footnotesize  Kog\u{a}lniceanu 1, 400084 Cluj-Napoca, Romania\\[0.1cm]
\footnotesize ${}^4$ Simion Stoilow Institute of Mathematics of the Romanian Academy,\\
\footnotesize Research group of the project PD-3-0152,\\
\footnotesize P. O. Box 1-764, RO-014700 Bucharest, Romania\\[0.1cm]
\footnotesize E-mails:  Laurentiu.Leustean@imar.ro, anicolae@math.ubbcluj.ro
}

\date{}

\maketitle

\begin{abstract}

In this note we point out that results on the asymptotic behaviour of an alternative iterative method are corollaries of corresponding results on the well-known Halpern iteration.  \\

\noindent {\em MSC:} 47J25, 47H09.
\end{abstract}

\maketitle

Let $X$ be a geodesic space, $C\subseteq X$ a convex subset and $T:C\to C$ a nonexpansive mapping. We consider the following iteration
\beq 
x_0\in C, \quad x_{n+1} = T(\lambda_{n+1}u + (1-\lambda_{n+1})x_n), \label{alternative-iteration}
\eeq
where $u\in C$ and $(\lambda_n)$ is a sequence $[0,1]$. This alternative iterative method was introduced for Banach spaces in \cite{Xu10} as a discretization of an approximating curve considered in \cite{ComHir06}. 
A more general version of this iteration is studied in \cite{ColLeuLopMar11}, while a  cyclic version of $(x_n)$ associated to a finite family of nonexpansive mappings is considered  in \cite{LeuNic14}.

Define for $n\geq 0$,
\beq
y_{n+1} = \lambda_{n+1}u + (1-\lambda_{n+1})x_n.
\eeq
Then for all $n\geq 1$, $x_n=Ty_n$. 
One can easily see by induction that the following holds.

\bprop
For all $n\geq 1$, $y_{n+1} = \lambda_{n+1}u + (1-\lambda_{n+1})Ty_n$.
\eprop
\details{$n=1$: $y_2=\lambda_2u + (1-\lambda_2)x_1=\lambda_2u +  (1-\lambda_2)T(\lambda_1u + (1-\lambda_1)x_0)=\lambda_2u + (1-\lambda_2)Ty_1$.

$n\Ra n+1$: \bua
y_{n+2} &=& \lambda_{n+2}u + (1-\lambda_{n+2})x_{n+1}=\lambda_{n+2}u + (1-\lambda_{n+2})T(\lambda_{n+1}u + (1-\lambda_{n+1})x_n)\\
&=&\lambda_{n+2}u + (1-\lambda_{n+2})Ty_{n+1}.
\eua
}
Thus, $(y_n)$ is the  well-known Halpern iteration \cite{Hal67,Wit92}. 
Since $T$ is nonexpansive, we have that for all $m,n\geq 1$ and for any fixed point $p$ of $T$,
\[d(x_m,x_n)\leq d(y_m,y_n) \quad \text{ and } \quad  d(x_n,p)\leq d(y_n,p).\]
Hence, one can obtain results concerning this iteration as immediate corollaries of corresponding results  on the Halpern iteration. 

\mbox{}

A first example is the following: if $(y_n)$ converges strongly to a fixed point of $T$, then $(x_n)$ also converges strongly to $p$.
It follows that \cite[Theorem 4.1]{Xu10}, the main result on the iteration $(x_n)$,  is a corollary of \cite[Theorem 3.1]{Xu02} and \cite{ShiTak97,Rei80}. Moreover, strong convergence results on the Halpern iteration, as well as their quantitative versions (see, e.g., \cite{Koh11,KohLeu12,Kor15,LeuNic15}), hold for the alternative iteration $(x_n)$ too.

\mbox{}

Let us give a second example. We say that a sequence $(z_n)$ is asymptotically regular if $\limn d(z_n,z_{n+1})=0$. A mapping $\Phi:(0,\infty)\to\N$ is a rate of asymptotic regularity for $(z_n)$ if for all $\eps>0$ and for all $n\geq \Phi(\eps)$, $d(z_n,z_{n+1})<\eps$.  
One has the following:
if $(y_n)$ is asymptotically regular with rate $\Phi$, then $(x_n)$ is asymptotically regular with the same rate.
As a consequence, (quantitative) asymptotic regularity results for the Halpern iteration hold without changes for the alternative iteration $(x_n)$. 

\mbox{ } 

\noindent
{\bf Acknowledgements:} \\[1mm] 
Lauren\c tiu Leu\c stean was supported by a grant of the Romanian 
National Authority for Scientific Research, CNCS - UEFISCDI, project 
number PN-II-ID-PCE-2011-3-0383. \\[1mm]
Adriana Nicolae was supported by a grant of the Romanian
Ministry of Education, CNCS - UEFISCDI, project number PN-II-RU-PD-2012-3-0152.

\end{document}